\documentclass[preprint,12pt]{elsarticle}

\usepackage{amssymb}
\usepackage{mathrsfs}
\usepackage{amsmath}
\usepackage{amsthm,amsfonts}
\usepackage{graphicx}
\usepackage{float}

\newtheorem{lemma}{Lemma}[section]

\newtheorem{theorem}[lemma]{Theorem}

\newenvironment{prof}[1][Proof]{\noindent\textit{#1}\quad }

\begin{document}

\newcommand{\js}{\hfill $\square$}

\begin{frontmatter}



\title{Tree-colorable maximal planar graphs
\footnote{This research is supported by  National Natural
Science Foundation of China (60974112, 30970960).}}




\author{ Enqiang Zhu$^{a,b}$, Zepeng Li$^{a,b}$, Zehui Shao$^{c,d}$, Jin Xu$^{a,b}$}

\address[PKU1]{School of Electronics Engineering and Computer Science, Peking University, Beijing, 100871, China}

\address[PKU2]{ Key Laboratory of High Confidence Software Technologies (Peking
University), Ministry of Education, CHINA}

\address[CDU1]{
  Key Laboratory of Pattern Recognition and Intelligent Information Processing, Institutions of Higher Education of Sichuan Province, China}

\address[CDU2]{
  School of Information Science and Technology, Chengdu University, Chengdu, 610106,  China\\
  }

\address{Peking University; Key Laboratory of High Confidence Software Technologies (Peking
University), Ministry of Education, CHINA\\
$Email: zhuenqiang@163.com$}

\begin{abstract}
A tree-coloring of a maximal planar graph is a proper vertex $4$-coloring such that every bichromatic subgraph, induced by this coloring, is a tree. A maximal planar graph $G$ is  tree-colorable if $G$ has a tree-coloring.
In this article, we prove that  a tree-colorable maximal planar graph $G$ with $\delta(G)\geq 4$ contains at least four odd-vertices.
 Moreover, for a tree-colorable maximal planar graph of minimum degree 4 that contains  exactly four odd-vertices, we show that the subgraph induced by its four odd-vertices is not a claw  and contains no triangles.
\end{abstract}

\begin{keyword}
Maximal planar graphs, tree-colorable maximal planar graphs, tree-coloring, claw, triangles.



\MSC 05C15
\end{keyword}

\end{frontmatter}


\section{Introduction}

The acyclic colorings was first studied by Gr\"{u}nbaum \cite{Grunbaum1973}, who wrote a long paper to research on the acyclic colorings of planar graphs. He proved that every planar graph is acyclic 9-colorable, and conjectured five colors are sufficient.
 Sure enough, three years later, Borodin \cite{Borodin1976}(also see \cite{Borodin1979}) gave a proof of Gr\"{u}nbaum's conjecture by showing that every planar graph is acyclic 5-colorable. In fact, this bound is the best possible for there exist planar graphs with no acyclic 4-colorings\cite{Grunbaum1973}. In 1973, Wegner \cite{Wegner1973} constructed a 4-colorable planar graph $G$, each 4-coloring of which possesses a cycle in every bichromatic subgraph. Afterwards Kostochka and Melnikov \cite{Kostochka1976}, in 1976, showed that graphs with no acyclic 4-coloring can be found among 3-degenerated bipartite planar graphs.

The research on acyclic $4$-colorable planar graphs always aroused  more attention. 
Some sufficient conditions have been obtained for a planar graph to be acyclic 4-colorable.
In 1999, Borodin, Kostochka, and Woodall \cite{Borodin1999} showed that planar graphs under the absence of 3- and 4-cycles are  acyclic 4-colorable; In 2006, Montassier, Raspaud, and Wang \cite{Montassier2006} proved that planar graphs, without 4-,5-, and 6-cycles, or without 4-, 5-, and 7-cycles, or without 4-, 5-, and intersecting 3-cycles, are  acyclic 4-colorable; In 2009, Chen and Raspaud \cite{Chen2009} proved that if a planar graph $G$ has no 4-, 5-, and 8-cycles, then $G$ is acyclic 4-colorable; Also in 2009, Borodin\cite{Borodin2009} showed that planar graphs without 4- and 6-cycles  are acyclic 4-colorable; Additionally, Borodin in 2011\cite{Borodin2011} and 2013\cite{Borodin2013} proved that planar graphs without 4- and 5-cycles are acyclic 4-colorable and acyclically 4-choosable, respectively.

\section{Preliminaries}

All of the graphs considered  are simple and finite.
For a graph $G$, we denote by $V(G)$, $E(G)$, $\delta(G)$ and $\Delta(G)$  the set of vertices, the set of edges, the \emph{minimum degree} and \emph{maximum degree} of $G$, respectively.
 For a vertex $u$ of $G$, $d_G(u)$ is the degree of $u$ in $G$. We call $u$ a \emph{$k$-vertex} if $d_G(u)=k$.
  If $k$ is an odd number, we say $u$ to be an \emph{odd-vertex}, and otherwise an \emph{even-vertex}.
  If $d_G(u)>0$, then each adjacent vertex of $u$ is called a \emph{neighbor} of $u$. The set of all neighbors of $u$ in $G$ is denoted by $N_G(u)$. Notice that $N_G(u)$ does not include $u$ itself. We then write $N_G[u]=N_G(u)\cup \{u\}$. For a subset $V'\subseteq V(G)$, denote by $G[V']$ the subgraph of $G$ induced by $V'$.
For more notations and terminologies, we refer the reader to the book \cite{Bondy2008}.

  A planar graph $G$ is called a \emph{plane triangulation} if the addition of any edge to $G$ results in a nonplanar graph. In this paper, triangulations are also known as \emph{maximal planar graphs}.

A \emph{k-coloring} of $G$ is an assignment of $k$ colors to $V(G)$ such that  no two adjacent vertices are assigned the same color. Alternatively, a $k$-coloring can be viewed as a partition $\{V_1,V_2,\cdots,V_k\}$ of $V$, where $V_i$ denotes the (possibly empty) set of vertices assigned color $i$, and is called a \emph{color class} of the coloring.

Let $f$ be a  coloring of a graph $G$, and $H$ be a subgraph of $G$.
We denote by $f(H)$ the set of colors assigned to $V(H)$ under $f$.
For a cycle $C$ of $G$, if $|f(C)|=2$, then we call $C$ a \emph{bichromatic cycle} of $f$, or say  $f$ \emph{contains bichromatic cycle} $C$. An \emph{acyclic $k$-coloring} of a graph $G$ is a \emph{k-coloring} with no bichromatic cycles  \cite{Grunbaum1973}.

For a maximal planar graph $G$, if $G$ has an acyclic 4-coloring $f$, then not only $f$ contains no bichromatic cycles,
but also any subgraph induced by  two color classes of $f$ is a tree.
So, it is more preferable to refer to such an acyclic 4-coloring as a \emph{tree-coloring} of $G$. Furthermore, if a maximal planar graph possesses a tree-coloring, then we say this graph is \emph{tree-colorable}.

The \emph{dual} graph $G^*$ of a plane graph $G$ is a graph that has a vertex corresponding to each face of $G$, and an edge joining two neighboring faces for each edge in $G$. It is well-known that the dual graphs of maximal planar graphs are planar cubic 3-connected graphs. Note that $G$ is a tree-colorable maximal planar graph if and only if its dual graph $G^*$ contains three Hamilton cycles such that each edges of $G^*$ is just contained in two of them. Since the problem of deciding whether a  planar cubic 3-connected graph contains a Hamilton cycle is NP-complete \cite{Garey1976}, we can deduce that the problem of deciding whether a maximal planar graph is tree-colorable is NP-complete. In addition, with regard to acyclic 4-colorability of planar graphs, it has been shown that acyclic 4-colorability is NP-complete for planar graphs with maximum degree 5,6,7, and 8 respectively and for planar bipartite graphs with the maximum degree 8 \cite{Mondala2013} \cite{Mondala2012} \cite{Ochem2005}.

As far as we know, there are no papers that have been written to study the tree-colorability (acyclic 4-colorability) of maximal planar graphs. Because maximal planar graphs contain a large number of 3-, 4-, or 5-cycles, we have reasons to believe that there exist lots of maximal planar graphs without tree-colorings. However, what are the characteristics of a tree-colorable maximal planar graph?
In this article,
we prove that  a tree-colorable maximal planar graph $G$ with $\delta(G)\geq 4$ contains at least four odd-vertices.
Furthermore, for a tree-colorable maximal planar graph of minimum degree 4 that contains  exactly four odd-vertices, we show that the subgraph induced by its four odd-vertices is not a claw  and contains no triangles.


\section{Main results}

First, we introduce a novel technique, named \emph{operation of contracting $4$-wheel}, which is very useful to the proof of the results throughout this paper.

A $\ell$-cycle $C$ is a cycle of length $\ell$. If $\ell$ is even, we call $C$ an even cycle, otherwise, an odd cycle.  A $n$-\emph{wheel $W_n$} (or simply \emph{wheel W}) is a graph with $n+1$ vertices ($n\geq 3$), formed by connecting a single vertex (called the \emph{center} of $W_n$) to all vertices of an $n$-cycle.

For a  maximal planar graph $G$ with $\delta(G)\geq 4$, it is obvious that any subgraph induced by a vertex and all of its neighbors is a wheel graph. Let $W$ be a 4-wheel subgraph of $G$.
The \emph{operation of contracting \emph{4}-wheel $W$} on $u,w$ of $G$, denoted by $\mathscr{D}_W^{u,w}(G)$, is to delete  $v$ from $G$ and identify vertices $u$ and $w$ (replace $u$, $w$ by a single vertex ($u,w$) incident to all the edges which were incident in $G$ to either $u$ or $w$), where $v$ is the center of $W$ and $u,w$ are two nonadjacent neighbors of $v$.
We denote by $\zeta_W^{u,w}(G)$ the resulting graph by conducting operation $\mathscr{D}_W^{u,w}(G)$. Clearly,
\begin{equation}
\begin{split}
 d_{\zeta_W^{u,w}(G)}((u,w))&=d_G(u)+d_G(w)-4, \\
 d_{\zeta_W^{u,w}(G)}(x)&=d_G(x)-2, \\
 d_{\zeta_W^{u,w}(G)}(y)&=d_G(y)-2,
\end{split}
\end{equation}
where $\{x,y\}=N_{G}(v)\setminus\{u,w\}$.
Notice that $\zeta_W^{u,w}(G)$ is still a maximal planar graph when $d_G(x)\geq 5$ and $d_G(y)\geq 5$.

\vspace{0.2cm}
We start with a few simple and useful conclusions.

\begin{lemma}\label{lemma1}
 Let $G$ be a tree-colorable maximal planar graph with a $4$-vertex $v$. Suppose that $f$ is a tree-coloring of $G$. Then $|f(N_G(v))|=3$, and $d_G(v_1)\geq 5$, $d_G(v_3)\geq 5$, where $v_1,v_3$ are the two nonadjacent neighbors of $v$ with $f(v_1)\neq f(v_3)$.
\end{lemma}
\begin{prof} Let $v_1,v_2,v_3,v_4$ be the four consecutive neighbors of $v$ in cyclic order.
It naturally follows that $|f(\{v_1,v_2,v_3,v_4\})|=3$ for $f$ is a tree-coloring.
Since $f(v_1)\neq f(v_3)$, we have $f(v_2)=f(v_4)$ and
$d_G(v_1)\geq 4$, $d_G(v_3)\geq 4$. If one of $v_1, v_3$ is a 4-vertex, say $v_1$, then it is unavoidable that $f$ contains a bichromatic cycle $v_2vv_4wv_2$ or $v_2v_3v_4wv_2$, where $\{w\}=N_G(v_1)\setminus \{v_2,v,v_4\}$. So $d_G(v_1)\geq 5$ and $d_G(v_3)\geq 5$. \qed
\end{prof}

\begin{lemma}\label{lemma2}
Let $G$ be a tree-colorable maximal planar graph with a $4$-vertex $v$, and $f$ be a tree-coloring of $G$. Then $\zeta_W^{v_1,v_2}(G)$ is still a tree-colorable maximal planar graph, where $W=G[ N_{G}[v]]$ and $v_1,v_2$ are two nonadjacent neighbors of $v$ such that $f(v_1)=f(v_2)$.
\end{lemma}
\begin{prof} By Lemma \ref{lemma1} $\delta(\zeta_W^{v_1,v_2}(G))\geq 3$ which implies  $\zeta_W^{v_1,v_2}(G)$ is still a maximal planar graph.
For any $v\in V(\zeta_W^{v_1,v_2}(G))$, if $v\neq (v_1,v_2)$, let $f^*(v)=f(v)$; otherwise, let $f^*(v)=f(v_1)$. Then, $f^*$ is a tree-coloring of  $\zeta_W^{v_1,v_2}(G)$.
\qed
\end{prof}

In this paper, we refer to the tree-coloring $f^*$ of $\zeta_W^{v_1,v_2}(G)$
in Lemma \ref{lemma2} as the \emph{inherited} tree-coloring of $f$.
Similar to the result of Lemma \ref{lemma2}, if a tree-colorable maximal planar graph $G$ contains 3-vertices, then the subgraph of $G$ obtained by deleting some (or all) 3-vertices is still a tree-colorable maximal planar graph.

Let $G$ be a graph with a cycle $C$.
We denote by $Int(C)$  the subgraph induced by $V(C)$ and all the vertices in the interior of $C$,  and denote by $Ext(C)$ the subgraph induced by $V(C)$ and vertices in the exterior of $C$.

A $k$-cycle $C$ of a connected graph $G$ is called a separating $k$-cycle if the deletion of $C$ from $G$ results in a disconnected graph.

\begin{lemma}\label{lemma3} A $3$-connected maximal planar graph $G$ is tree-colorable if and only if  for any separating $3$-cycle  $C$ of $G$, both of $Int(C)$ and $Ext(C)$ are tree-colorable.
\end{lemma}
\begin{prof}
This result is obvious, so we omit the proof. \qed
\end{prof}

\vspace{0.2cm}

Based on the above tree lemmas, we  give the first main result of this section as follow.
\begin{theorem}\label{theorem 2.4} A tree-colorable maximal planar graph of minimum degree at least $4$ contains at least four odd-vertices.
\end{theorem}
\begin{prof}
Let $G$ be a tree-colorable maximal planar graph with $\delta(G)\geq 4$. Then the minimum degree of $G$ is either 4 or 5. Indeed, it suffices to consider the case of $\delta(G)=4$ because $G$ contains at least twelve 5-vertices by the Euler Formula when $\delta(G)=5$.

If the conclusion fails to hold when $\delta(G)=4$, let $G^\prime$ be a counterexample on the fewest vertices to the theorem, i.e. $G^\prime$ is a tree-colorable maximal planar graph of $\delta(G^\prime)=4$ with $o(G^\prime)<4$,  where $o(G^\prime)$ is the number of odd-vertices of $G^\prime$. It is obvious that $o(G^\prime)=2$ or $o(G^\prime)=0$. Thus, by using the well-known relation
\[
\sum\limits_{v\in V(G)}(d(v)-6)=-12,
\]
we can deduce $G^\prime$ contains at least five 4-vertices.

Let $f$ be an arbitrary tree-coloring of  $G^\prime$. If $G^\prime$ contains no 5-vertices, then for any 4-vertex $u$ and its two nonadjacent neighbors $u_1$,$u_2$ with $f(u_1)=f(u_2)$,  $\zeta_W^{u_1,u_2}(G^\prime)$ is still a tree-colorable maximal planar of minimum degree at least 4 and contains at most two odd-vertices by formula (1), where $W=G^\prime[N_{G^\prime}[u]]$. This contradicts the assumption of $G^\prime$. So we only need to consider the case that $G^\prime$ contains 5-vertices.

Note that for any 5-vertex $v$ of $G^\prime$, there are at most three 4-vertices in $N_{G^\prime}(v)$. Otherwise, if there are four (or five) 4-vertices in $N_{G^\prime}(v)$, then $G^\prime$  is the graph $G_7$ shown in Figure \ref{figure1}(a). However,  it is an easy task to prove that $G_7$ contains no tree-colorings, and a contradiction. We now turn to show that there are also no three vertices in $N_{G^\prime}(v)$ with degree 4. If not, let $v_1,v_2,v_3$ be three 4-vertices of $N_{G^\prime}(v)$.

\begin{figure}[H]
  \centering
  \includegraphics[width=300pt]{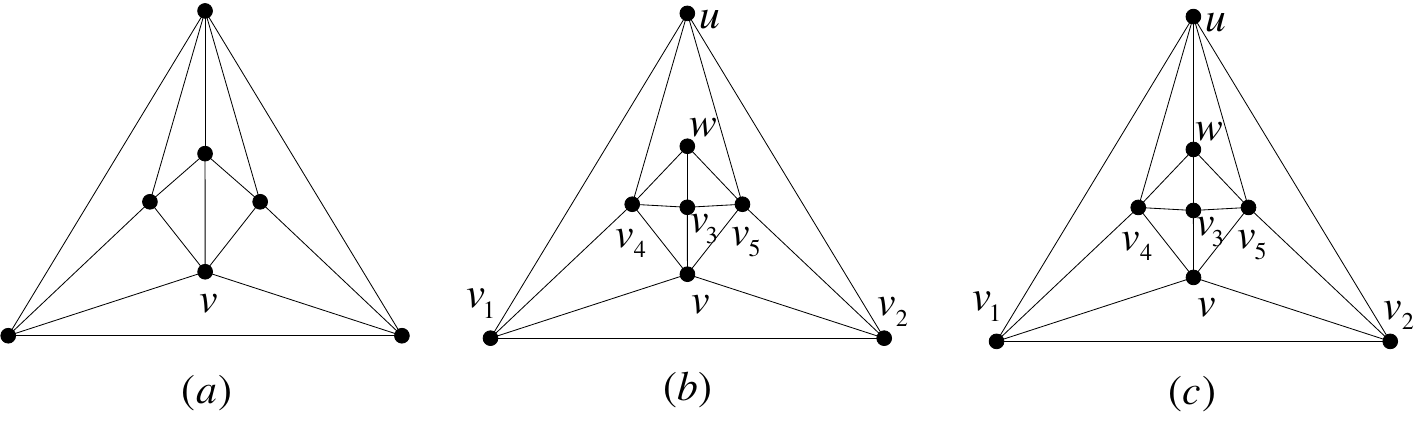}\\
\caption{$(a) G_7, (b) H, (c) G_8$}\label{figure1}
\end{figure}

(1) $v_1,v_2,v_3$ are three consecutive vertices, i.e. $G^\prime[\{v_1,v_2,v_3\}]$ contains two edges. However, it is readily to check that $G^\prime$ contains subgraph $G_7$, which contradicts the assumption that $G^\prime$ contains tree-coloring.

(2) $G^\prime[\{v_1,v_2,v_3\}]$ contains only one edge, w.l.o.g. say $v_1v_2\in E(G^\prime)$. Then $G^\prime$ contains a subgraph $H$ isomorphic to the graph shown in Figure \ref{figure1}(b). It is easy to see that $f(v)=f(w)$.

If $d_G^\prime(v_4)\geq 6$ and $d_G^\prime(v_5)\geq 6$, then $\zeta_W^{v,w}(G^\prime)$ is still a tree-colorable maximal planar graph of minimum degree at least 4 and contains at most two odd-vertices by Lemma \ref{lemma2}, where $W=G^\prime[N_{G^\prime}[v]]$, and a contraction with minimum property of $G^\prime$.

If there is a 5-vertex in $\{v_4,v_5\}$, say $v_5$, then $G^\prime$ is either the graph $G_8$ shown in Figure \ref{figure1}(c) that contains four 5-degree vertices (a contradiction with $G^\prime$), or a 3-connected graph with separating 3-cycle $C=v_4wuv_4$. For the latter case, either $\delta(Int(C))\geq 4$, or there exists another separating 3-cycle $C'$ in $Int(C)$ such that $\delta(Int(C'))\geq 4$ (because there must be a separating 3-cycle $C'$ such that $Int(C')$ is 4-connected, otherwise there are 3-vertices in $Int(C')$). By Lemma \ref{lemma3}, $Int(C)$ ( or $Int(C')$) is a tree-colorable maximal planar graph with minimum degree at least 4 and contains at most two odd-vertices, which contradicts the assumption of $G^\prime$.

The above two cases imply that any $5$-vertex $v$ in $G^\prime$ has  at most two neighbors with degree 4. Since there are at least five 4-vertices in $G^\prime$, we can always find a $4$-vertex $v'$ such that $N_{G^\prime}(v')$ contains no 5-vertices. So, by Lemma \ref{lemma1} and \ref{lemma2}, $\zeta_W^{v'_1,v'_2}(G^\prime)$ is still a tree-colorable maximal planar graph of minimum degree at least 4, where $W=G^\prime[N_{G^\prime}[v']]$ and $v'_1,v'_2\in N_{G^\prime}(v')$ with $f(v'_1)=f(v'_2)$. However, $\zeta_W^{v'_1,v'_2}(G^\prime)$ contains at most two odd-vertices, and this contradicts the choice of $G^\prime$.\qed
\end{prof}

By Lemma \ref{lemma3}, it clearly suffices to consider tree-colorable  maximal planar graphs without separating 3-cycle.
In what follows, we denote by $MPG4$ the class of tree-colorable 4-connected maximal planar graphs with exact four odd-vertices. Furthermore, for a graph $G\in MPG4$, we denote by $V^4(G)$ the set of the four odd-vertices of $G$. Obviously, the minimum degree of graphs in $MPG4$ is 4.  Now, we turn to discuss the structural properties of graphs in $MPG4$.

For a graph $G$ in $MPG4$ and a 4-vertex $v$, if there are two
vertices $v_1, v_2\in N_G(v)$ such that $v_1v_2 \notin E(G)$
and $\zeta_W^{v_1,v_2}(G)$ 
is still a graph in $MPG4$, then we refer to such vertex $v$ as a \emph{contractible vertex} of $G$.

In order to investigate the structure of the subgraph induced by the four odd-vertices of a graph in $MPG4$, we need a lemma as follow.

\begin{lemma}\label{lemma4}
Let $G$ be a graph in $MPG4$.

\begin{description}
  \item[(1)] If $G$ contains a $5$-vertex $v$ such that $N_G(v)$ contains at least three $4$-vertices, then either $G$ is the graph isomorphic to $G_7$ or $G_8$, or $G$ contains contractible vertices.
  \item[(2)] If $G$ contains a $7$-vertex $v$ such that $N_G(v)$ contains at least five $4$-vertices, then $G$ contains contractible vertices.
  \item[(3)] If $G$ contains a $9$-vertex $v$ such that $N_G(v)$ contains at least six $4$-vertices, then either $G$ has contractible vertices, or $G$ is the graph isomorphic to Figure $\ref{figure2.5}$.
\end{description}

\end{lemma}
\begin{prof}
\textbf{(1)}. According to the proof of Theorem \ref{theorem 2.4}, we can know that  $G$ contains either subgraph $G_7$ or subgraph $H$. Since $G$ is 4-connected, it follows  that either $G$ is the graph isomorphic to $G_7$ or $G_8$, or $G$ contains contractible vertices (see the vertex $v_3$ of graph $H$ shown in Figure \ref{figure1}(b)).

\textbf{(2)}. Let $v$ be a 7-vertex of $G$, and $N_G(v)=\{v_1,v_2,v_3,v_4,v_5,v_6,v_7\}$ in cyclic order. If $N_G(v)$ contains at least five $4$-vertices, then at least
three of them are consecutive, say $v_1,v_2,v_3$. Denote by $v_8$ the common neighbour (except $v$) of them, and then we have $d_G(v_8)\geq 6$. Otherwise, $G$ contains separating 3-cycle $v_7vv_4v_7$. By Lemma \ref{lemma1} for each tree-coloring $f$ of $G$, we have $f(v_1)=f(v_3)$. So $v_2$ is a contractible vertex.

\textbf{(3)}. Let $v$ be a 9-vertex of $G$, and $N_G(v)=\{v_1,v_2,v_3,v_4,v_5,v_6,v_7,v_8,v_9\}$ in cyclic order (see Figure \ref{figure2.5}). If there are three consecutive 4-vertices in  $N_G(v)$, similarly to \textbf{(2)} $G$ has contractive vertices. If there are no three consecutive 4-vertices in $N_G(v)$, then the number of 4-vertices of $N_G(v)$ is exactly 6. W.o.l.g. we assume $v_2,v_3,v_5,v_6,v_8,v_9$ are the six 4-vertices. Because $G$ is 4-connected, we can assume that the common neighbor(except $v$) of $v_2$ and $v_3$ is $u_1$, the common neighbor (except $v$) of $v_5$ and $v_6$  is $u_2$, and  the common neighbor (except $v$) of $v_8$ and $v_9$  is $u_3$ (see Figure \ref{figure2.5}). If one of $u_1,u_2,u_3$ is a 6-vertices, say $u_1$,  then $v_2$ and $v_3$ are contractible vertices. If $d_G(u_1)=d_G(u_2)=d_G(u_3)=5$, then it follows that $G$ is the graph isomorphic to Figure $\ref{figure2.5}$. \qed

\begin{figure}[H]
  \centering
  \includegraphics[width=120pt]{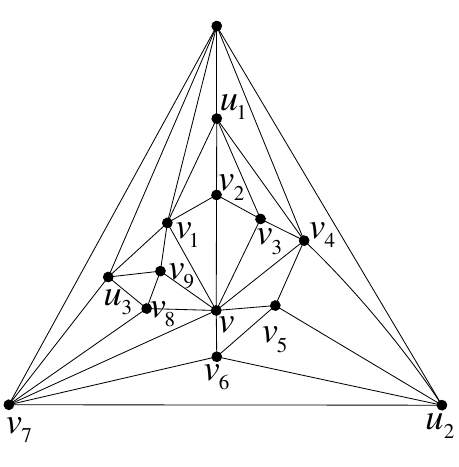}\\
  \caption{A graph}\label{figure2.5}
\end{figure}

\end{prof}

We then prove that the subgraph induced by the four odd-vertices  of a graph in $MPG4$ contains no triangles.

\begin{theorem}\label{theorem 2.5} Let $G$ be a  graph in $MPG4$ with $n$ vertices. Then  $G[V^4(G)]$ contains no triangles.
\end{theorem}

\begin{prof}
With the help of the software \emph{plantri} developed by \emph{McKay} \cite{McKay2007}, we confirm that
there are 1,0,2 and 1 graphs in $MPG4$ when $n=8,9,10$ and $11$, respectively.
We now proceed by induction on $n$.

Suppose that the theorem holds for all graphs in $MPG4$ with fewer than $n$($\geq 12$) vertices.  Let $G$ be a graph  in $MPG4$ with $n$ vertices, and $V^4$=\{$u_1,u_2,u_3,u_4$\}.
We claim that $G[V^4(G)]$ contains no triangles.
If not, we w.l.o.g. assume $u_1u_2u_3$ is a triangle of $G[V^4(G)]$. Then $G$ contains no contractible vertices. Otherwise let $u$ be a contractible vertex, i.e. there exist two vertices $x_1,x_2$ in $N_G(u)$ such that $x_1x_2\not\in E(G)$ and $\zeta_{W'}^{x_1,x_2}(G) \in $ $MPG4$, where $W'=G[N_G[u]]$. However it is an easy task to show that the subgraph of  $\zeta_{W'}^{x_1,x_2}(G)$ induced by its four odd-vertices also contains a triangle, and this contradicts the hypothesis.

Notice that $5\leq d_G(u_4)\leq 9$. Otherwise, if $d_G(u_4)\geq 11$, then $G$ contains at least seven 4-vertices. This indicates that there exists a 4-vertex adjacent no 5-vertices  by Lemma \ref{lemma4} (1). So the 4-vertex is a contractible vertex.

If $d_G(u_4)=5$, then $G$ contains at least four 4-vertices, and $N_G(u_4)$ contains at most two 4-vertices by Lemma \ref{lemma4} (1); If $d_G(u_4)=7$, then $G$ contains at least five 4-vertices, and $N_G(u_4)$ contains at most four 4-vertices by Lemma \ref{lemma4} (2); If $d_G(u_4)=9$, then $G$ contains at least six 4-vertices and $N_G(u_4)$ contains at most four 4-vertices by Lemma \ref{lemma4} (3).  So,
  we can always find a 4-vertex, say $v'$, such that $u_4\not \in N_G(v')$. Let $v_1$,$v_2$,$v_3$,$v_4$ be  the four consecutive neighbors of $v'$ (see Figure \ref{figure2}(a)). We now assume $f(v_1)=f(v_3)$ for any tree-coloring $f$ of $G$, and then $f(v_2)\neq f(v_4)$ and $d_{G}(v_2)\geq 5$, $d_{G}(v_4)\geq 5$ by Lemma \ref{lemma1}. In terms of the relation between $\{u_1,u_2,u_3\}$ and $N_G(v')$, there are three cases which can happen. Obviously, $\{u_1,u_2,u_3\}\not \subset N_G(v')$.

\begin{figure}[H]
  \centering
  \includegraphics[width=380pt]{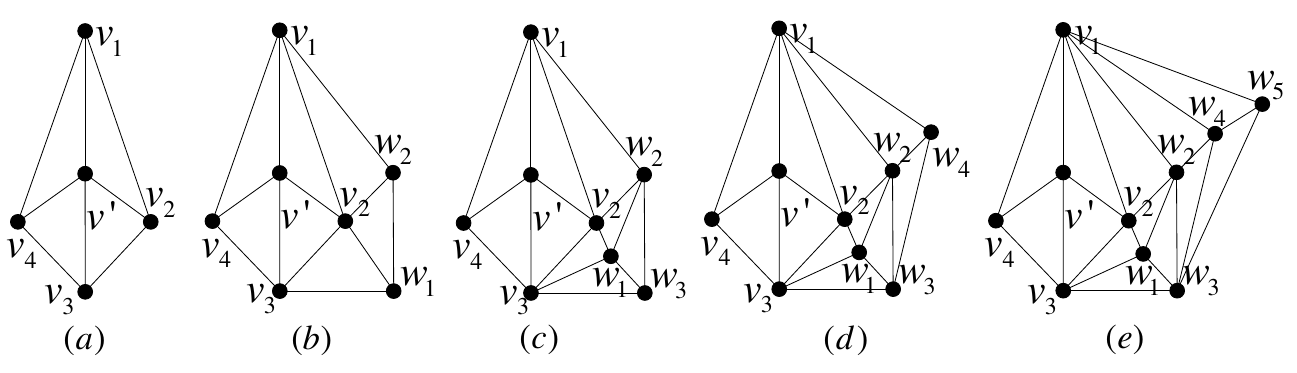}\\
  \caption{}\label{figure2}
\end{figure}

$Case ~1.$ One of  $\{u_1,u_2,u_3\}$ belongs to $N_G(v')$, say $u_1$.
By symmetry, it is sufficient to consider $u_1=v_1$ or  $u_1=v_2$.

If  $u_1=v_1$,  then $d_{G}(v_2)>5$ and $d_{G}(v_4)>5$. So $v'$ is a contractible vertex.

If $u_1=v_2$, it follows $d_G(u_1)=5$ (otherwise $v'$ is a contractible vertex of $G$).
Considering that $N_G(v_2=u_1)=\{v_1,v',v_3,u_3,u_2\}$, we have that $\zeta_{W_1}^{v_1,v_3}(G)-v_2$ is a tree-colorable maximal planar graph with minimum degree at least 4 with only two odd-vertices, where $W_1=G[N_G[v']]$. This contradicts Theorem \ref{theorem 2.4}.

$Case ~2.$ Two vertices of $\{u_1,u_2,u_3\}$ belong to $N_G(v')$, w.l.o.g. let $v_1=u_1,v_2=u_2$. If $d_G(v_2)\geq 7$, then $v'$ is a contractible vertex.
If $d_G(v_2)=5$, let $N_G(v_2)=\{u_1,v',v_3,w_1,w_2\}$, where $w_2=u_3$ (see Figure \ref{figure2}(b)).

$Case ~2.1.$ $d_{G}(w_1)\geq 5$, then $\zeta_{W_1}^{v_1,v_3}(G)-v_2$ is still a tree-colorable maximal planar graph with minimum degree 4, but contains at most two odd-vertices, and a contradiction with Theorem \ref{theorem 2.4}.

$Case ~2.2.$ $d_{G}(w_1)=4$, let $N_G(w_1)=\{v_3,v_2,w_2,w_3\}$ (see Figure \ref{figure2}(c)). Obviously, $d_{G}(v_3)\geq 6$ and $f(w_3)=f(v_2)$. If $d_{G}(w_2)\geq 7$, $w_1$ is a contractible vertex. If $d_{G}(w_2)=5$, let $N_G(w_2)=\{v_1,v_2,w_1,w_3,w_4\}$ (see Figure \ref{figure2}(d)). Then $\zeta_{W_1}^{v_2,w_3}(G)-w_2$ is tree-colorable maximal planar graph with minimum degree at least 4 when $d_{G}(w_4)\geq 5$, but contains at most two odd-vertices. This contradicts to Theorem \ref{theorem 2.4}; When $d_{G}(w_2)=5$ and $d_{G}(w_4)=4$, let $N_G(w_4)=\{v_1,w_2,w_3,w_5\}$ (see Figure \ref{figure2}(e). Noting that here $w_4v_4$ is not an edge of $G$. Otherwise, $w_3v_4$ is also an edge of $G$ for $d_{G}(w_4)=4$, and $G$ is a maximal planar graph of order 9 and contains more than four 5-degree vertices). Clearly, $d_{G}(v_1)\geq 7$ and $f(w_5)=f(w_2)=f(v_4)$ or $f(v')$, so $w_5v_3\not\in E(G)$ and $d_{G}(w_3)\geq 6$. This implies that $w_4$ is contractible vertex.

All of the above discussions show that $G[V^4(G)]$ contains no triangles.
 \qed
\end{prof}

Recall that a star $S_k(k\geq 2)$ is the complete bipartite graph $K_{1,k}$, which is a tree with one internal node and $k$ leaves. A star with $3$ edges is called a \emph{claw}, i.e. $S_3$. We now in a position to show that the  subgraph induced by the four odd-vertices  of a graph in $MPG4$ is not a claw.

\begin{lemma}\label{lemma2.6} Suppose that $G$ is a $4$-connected maximal planar graph satisfying the following three restrictions.
\begin{enumerate}[1)]
  \item Except one 9-vertex, three $5$-vertices, and six $4$-vertices, all of other vertices of $G$ are $6$-vertices;
  \item Any two 5-vertices are nonadjacent each other, and all $5$-vertices are  neighbors of the $9$-vertex, and every 5-vertex is adjacent to exactly two $4$-vertices;
  \item Each $4$-vertex is adjacent to one $5$-vertex.
\end{enumerate}
Then $G$ is a graph isomorphic to one of the graphs shown in Figure \emph{\ref{figure2.6}(a), (b), (c), (d),(e)}.

\begin{figure}[H]
  \centering
 \includegraphics[width=350pt]{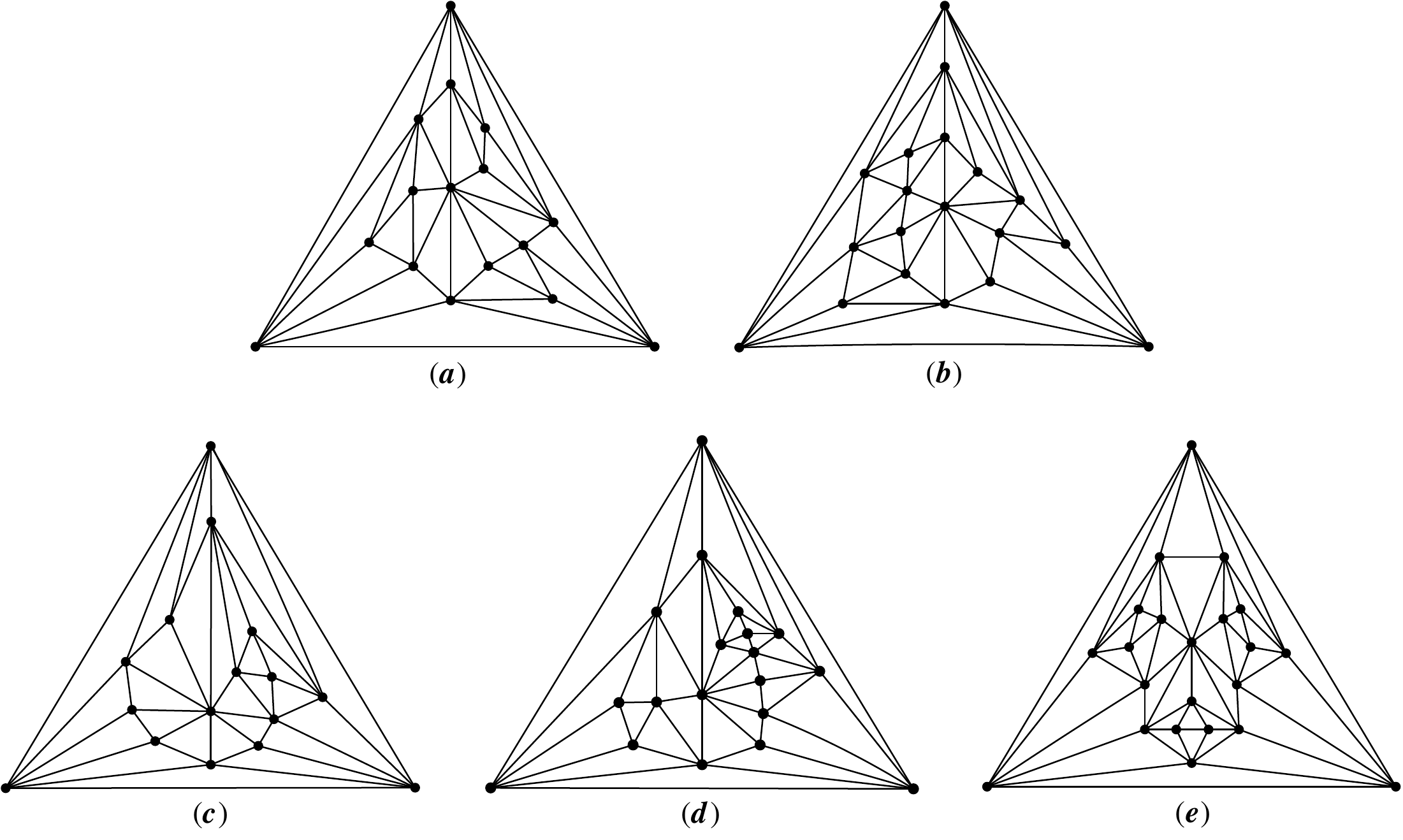}\\
 \caption{Five unavoidable graphs of Lemma\ref{lemma2.6}}\label{figure2.6}
\end{figure}

\end{lemma}

\begin{prof} Let $u_0$ be the 9-vertex.  Since $G$ is 4-connected, $G[N_G(u_0)]$ is a cycle $C$, denoted by $C=u_1u_2 \cdots u_9u_1$ (see Figure\label{figure2.6.1}(a)). We first show that there are no two 5-vertices has a common neighbor on $C$. If not, w.l.o.g.  we assume $u_2$, $u_9$ are two 5-vertices. Clearly, $u_1$ is a 6-vertex, and $u_2$, $u_9$ have no common neighbors (except $u_0$), see Figure \ref{figure2.6.1} (a), where $v_2,v_3$ (resp. $v_4,v_5$) are the neighbors of $u_2$ (resp. $u_9$) not on $C$, and $v_1$ is a neighbor of $u_1$ not on $C$. Obviously, $d_G(v_1)=6$. As each 5-vertex has exactly two neighbors of degree 4, we consider the following three cases.

$Case~1.$  $d_G(v_2)=d_G(v_3)=4$, i.e. $v_1v_3$, $v_1u_3\in E(G)$. Then it is impossible  $d_G(v_4)=d_G(v_5)=4$, otherwise, $d_G(v_1)\geq 7$.

$Case~1.1.$  $d_G(v_4)=d_G(u_8)=4$, i.e. $v_1v_5$, $v_5u_7\in E(G)$, see Figure \ref{figure2.6.1}(b). For $v_1$ is a 6-vertex, we have $v_5u_3\in E(G)$, which implies one of $u_3$ and $v_5$ is a vertex of degree at least 7, and a contradiction with $G$.

$Case~1.2.$  $d_G(v_5)=d_G(u_8)=4$, i.e. $v_4u_7$, $v_5u_7\in E(G)$, see Figure \ref{figure2.6.1}(c). For $u_7,u_3,v_4,v_1$ are 6-vertices, we can know $u_5$ is a 3-vertex under the condition $d_G(u_7)=d_G(u_3)=d_G(v_4)=d_G(v_1)=6$, and  a contraction with $G$.

\begin{figure}[H]
  \centering
 \includegraphics[width=340pt]{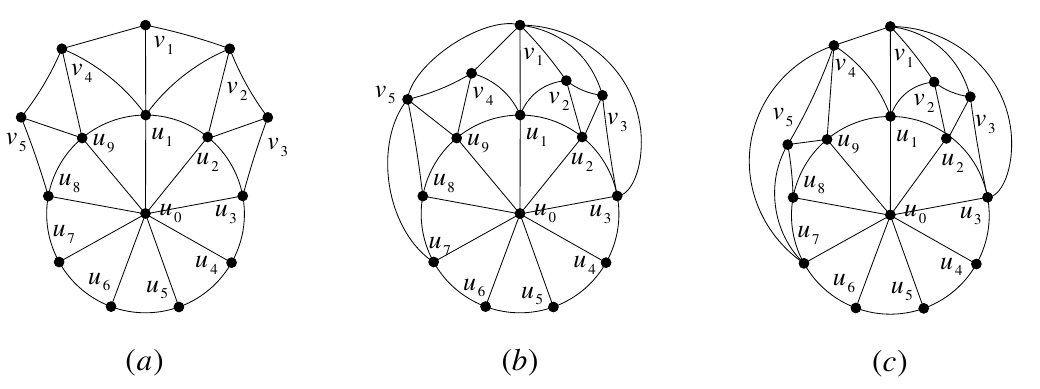}\\
 \caption{Graphs for Case 1}\label{figure2.6.1}
\end{figure}

$Case~2.$  $d_G(v_2)=d_G(u_3)=4$, see Figure \ref{figure2.6.2}(a).

$Case~2.1.$  $d_G(v_4)=d_G(v_5)=4$, then $v_1v_5$, $v_1u_8\in E(G)$, see Figure \ref{figure2.6.2}(a). For $v_1$ is a 6-vertex, at least one of $u_8$ and $v_3$ is a vertex of degree at least 7, and this contradicts to $G$.

$Case~2.2.$  $d_G(v_4)=d_G(u_8)=4$, then
$v_1v_5$, $v_5u_7\in E(G)$, see Figure \ref{figure2.6.2}(b). In this case, for $v_1,v_3,v_5$ are 6-vertices, and it can be seen that  there are no edges between $v_1$ and $u_4,u_5,u_6,u_7$ respectively ( By symmetry, we only consider $v_1u_4\not\in E(G)$ and $v_1u_5\not\in E(G)$. If $v_1u_4\in E(G)$, $d_G(v_3)=5$; if $v_1u_5\in E(G)$, $d_G(u_6)=3$ based on $d_G(v_1)=d_G(v_3)=6$. We  denote the additional neighbor of $v_1$ by $v_6$, see Figure \ref{figure2.6.2}(c). By $d_G(v_3)=d_G(v_5)=d_G(v_6)=6$ and $d_G(u_4)=d_G(u_7)=6$, we can further known  $d_G(u_5)=d_G(u_6)=4$, and a contraction with the condition that there are three 5-vertices on $C$.

$Case~2.3.$  $d_G(v_5)=d_G(u_8)=4$, then
$v_5u_7$, $v_4u_7\in E(G)$ and $d_G(u_7)=6$. Since $v_1,v_3,v_4$ are 6-vertices, we can know that $u_7v_1,u_7v_3\not\in E(G)$. Considering the additional neighbor of $u_7$, denote by $v_6$ (see Figure \ref{figure2.6.2}(d)). If $d_G(v_6)=4$, $u_6$ will be a 6-vertex since $v_1,v_3$ are 6-vertices, which implies $v_6$ is not adjacent to a 5-vertex and a contradiction. So  $v_6$ is a 6-vertex, i.e. $v_1u_6\not\in E(G)$. Further, as $v_3$ is a 6-vertex, there are no edges between $v_1$ and $u_4,u_5,u_6$. Let $v_7$ be additional neighbor of $v_3$. Because $v_6$, $u_4$, $v_7$ are 6-vertices, if $d_G(u_6)=5$ we have $d_G(u_5)=5$ and a contradiction with $G$; If $d_G(u_6)=6$ and $d_G(u_5)=5$, then $G$ is the graph isomorphic to the graph shown in Figure \ref{figure2.6.2}(d), and a contradiction with $G$.

\begin{figure}[H]
  \centering
 \includegraphics[width=400pt]{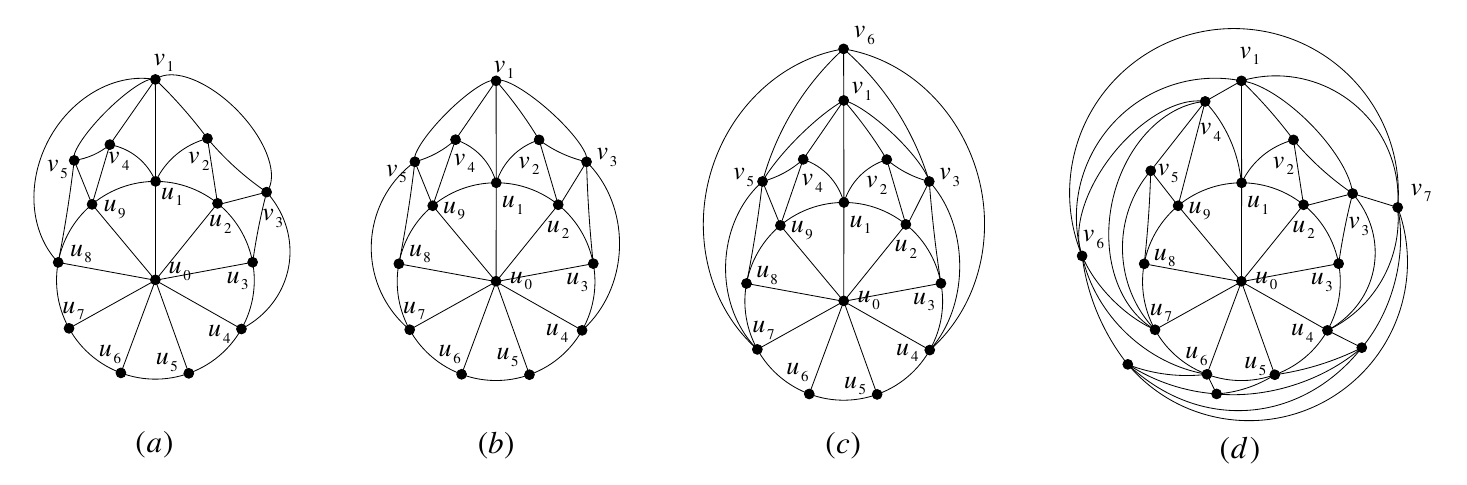}\\
 \caption{Graphs for Case 2}\label{figure2.6.2}
\end{figure}

$Case~3.$  $d_G(v_3)=d_G(u_3)=4$. By symmetry, we need only to check the case that $d_G(v_3)=d_G(u_3)=4$ and $d_G(v_5)=d_G(u_8)=4$. With the analogously arguing process, it is also to show that this case fails to exist.

Based on the above analysis, we confirm that the three 5-vertices of $G$ are  equably distributed on $C$, i.e. no two of them have a common neighbor on $C$. In what follows, w.l.o.g. we assume $u_1,u_4,u_7$ are the three 5-vertices of $G$ on $C$. If there are 4-vertex on $C$, suppose w.l.o.g. that $u_2$ is a 4-vertex. Let $v_1,v_2$ be additional two neighbors of $u_1$, where $v_1$ is the common neighbor of $u_1$ and $u_2$.

 (1) $d_G(v_1)=4$. Since $d_G(u_3)=6$ and $d_G(u_4)=5$, $u_4$ has additional two neighbors, say $v_3,v_4$, where $v_3$ is the common neighbor of $u_3$ and $u_4$ (see Figure \ref{figure2.6.3}(a)).
 Obviously, it is impossible that $d_G(v_3)=d_G(v_4)=4$, otherwise $v_2$ would be a vertex of degree at least 7.
        If $d_G(v_3)=d_G(u_5)=4$, then $d_G(u_7)=3$ for $v_2,v_4,u_6,u_9$ are 6-vertices; If $d_G(v_4)=d_G(u_5)=4$, then $G$ is the graph isomorphic to Figure \ref{figure2.6}(a).

  (2) $d_G(v_2)=4$. We claim that $v_2u_i\not\in E(G)$ for $i=3,4,5,6,7,8$. Since $d_G(v_1)=d_G(u_3)=6$, it is indirectly $v_2u_3,v_2u_4\not\in E(G)$. If $v_2u_5\in E(G)$, then $v_1u_5\in E(G)$ that indicates $u_5$ is a 6-vertex. Hence we have $d_G(u_3)=5$, a contradiction. Similarly, we have $v_2u_i\not\in E(G)$ for $i=6,7,8$.
         Let $v_3$ be another neighbor of $v_2$ and $v_4$ be the common neighbor of $v_3$ and $u_9$ (see Figure \ref{figure2.6.3}(b)). Also, there are no edges between $v_1$ and $v_4, u_4,u_5,u_6,u_7,u_8$ since $v_3,u_3$ are 6-vertices, and $u_i$ is a 6-vertex if $v_1u_i\in E(G)$ for $i=5,6,7,8$. So, there is another neighbor of $v_1$, say $v_5$ (see  Figure \ref{figure2.6.3}(b)). Since $u_3$ and $u_4$ are 6-vertex and 5-vertex respectively and $d_G(v_3)=d_G(v_4)=6$, $u_4$ has additional two neighbors, say $v_6,v_7$, where $v_6$ is the common neighbor of $u_3$ and $u_4$ (see Figure \ref{figure2.6.3}(b)).
         If $v_7,u_5$ are 4-vertices, then $u_7$ does not contain two neighbors of degree 4 based on $d_G(v_3)=d_G(v_5)=d_G(v_6)=6$; If $ v_6,v_7$ are 4-vertices, then one of $v_3,u_5$ is a vertex of degree at least 7 on the basis of $d_G(v_5)=6$; If $v_6,u_5$ are 4-vertices, then $G$ is isomorphic to the graph shown in Figure \ref{figure2.6}(b);

(3) $d_G(u_9)=4$, then $d_G(v_1)=d_G(v_2)=6$. Let $v_3,v_4$ be the other two neighbors of $u_3$, obviously $d_G(v_3)=d_G(v_4)=6$, where $v_3$ is the common neighbor of $v_1$ and $u_3$. So, $u_4$ is not adjacent to $v_1,v_2,v_3$ and $v_4$. Suppose the additional neighbor of $u_4$ is $v_5$, see Figure \ref{figure2.6.3}(c).  If $d_G(v_4)=d_G(u_5)=4$, then $u_7$ does not contain two neighbors of degree 4 based on $d_G(v_1)=d_G(v_2)=d_G(v_3)=d_G(v_5)=d_G(u_6)=6$; If $d_G(v_5)=d_G(u_5)=4$, then also $u_7$ does not contain two neighbors of degree 4 based on $d_G(v_1)=d_G(v_2)=d_G(v_3)=d_G(v_4)=6$; If $d_G(v_4)=d_G(v_5)=4$, then $G$ is a graph isomorphic to either the graph shown in Figure \ref{figure2.6}(c) or the graph shown in Figure\ref{figure2.6}(d).

 The above discussions show that when there are 4-vertices on $C$, $G$ is the graphs isomorphic to the graphs shown in Figure \ref{figure2.6}(a),(b),(c),(d). Moreover, if there are no 4-vertices on $C$, then it is obvious that $G$ is the graph isomorphic to the Figure \ref{figure2.6}(e).
\begin{figure}[H]
  \centering
 \includegraphics[width=360pt]{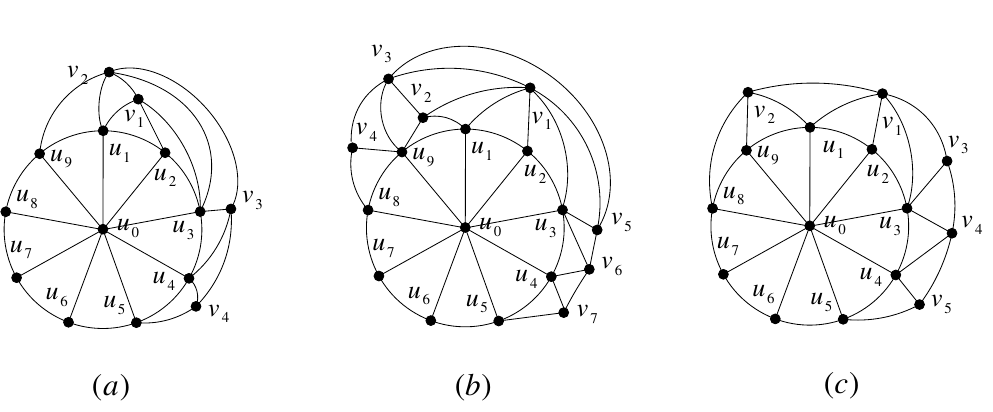}\\
 \caption{Graphs for Case 3}\label{figure2.6.3}
\end{figure}

\end{prof}

\begin{theorem}\label{theorem2.7}  Let $G$ be a graph in $MPG4$, and $V^4(G)=\{u_1,u_2,u_3,u_4\}$. Then $G[V^4(G)]$ is not a claw.
\end{theorem}

\begin{figure}[H]
  \centering
  \includegraphics[width=120pt]{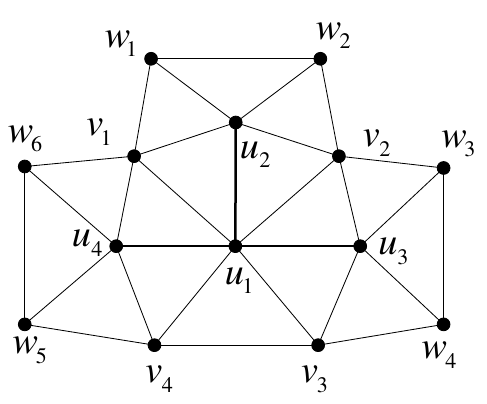}\\
  \caption{A subgraph}\label{figure3}
\end{figure}

\begin{prof}
If the result fails to hold, select a minimum counterexample $G'$, i.e. $G'$ is a graph in $MPG4$ with
the fewest vertices such that $G'[\{u_1,u_2,u_3,u_4\}]$ is a claw(see Figure \ref{figure3}(a)).
 Obviously, $G'$ does not contain contractible vertices.
 Thus,  it suffices to consider no cases that 5-vertex contains at least three neighbors with degree 4 by Lemma \ref{lemma4}.
Furthermore, according to Theorem \ref{theorem 2.5}, it follows $d_{G'}(u_1)\geq 7$.

If $d_{G'}(u_1)=7$, then $G'$ contains at least five 4-vertices. So $G'$ contains  a contractible vertex when there is at least one 7-vertex in  $\{u_2,u_3,u_4\}$; However, if $d_{G'}(u_2)=d_{G'}(u_3)=d_{G'}(u_4)=5$, (w.l.o.g. see Figure \ref{figure3}), then $d_{G'}(v_1)\geq 6$ and $d_{G'}(v_2)\geq 6$. Indeed, if $V'=\{w_1,w_2,w_3,w_4,w_5,w_6\,v_3,v_4\}$ contains $\ell$(=5,6,7) 4-vertices, then there are at least two vertices $x,y$ in $N_G[V']$ such that $d_G(x)+d_G(y)\geq 2\ell +4$ by the properties of $G'$ (All of vertices of $G'$ except $u_1,u_2,u_3,u_4$ are even-vertices) in this case.
So there are at least one 4-vertex without neighbors of degree 5, which means that $G'$ contains a contractible vertex and this contradicts the choice of $G'$.

If $d_{G'}(u_1)=9$, then there are at least six 4-vertices in $G'$. In this case, because each 5-vertex contains at most two neighbors of degree 4, it suffices to consider the unique case: $G'$ contains exact six 4-vertices, $d_{G'}(u_2)=d_{G'}(u_3)=d_{G'}(u_4)=5$, and all other vertices of $G'$ have degree 6. Otherwise, $G'$ contains contractible vertices. Thus, it requires that $N_G(\{u_2,u_3,u_4\})$ contains six 4-vertices and
each of $u_2,u_3,u_4$ contains exactly two distinct 4-vertices. Then, $G'$ is one of the graphs shown in Figure \ref{figure2.6} by Lemma \ref{lemma2.6}, which is not a tree-colorable maximal planar graph.

If $d_{G'}(u_1)\geq 11$, then $G'$ contains at least seven 4-vertices. So at least one 4-vertex has no neighbors of degree 5, and this 4-vertex is a contractible vertex of $G'$, and a contradiction. \qed
\end{prof}

\vspace{0.2cm}

\begin{figure}[H]
  \centering
  \includegraphics[width=300pt]{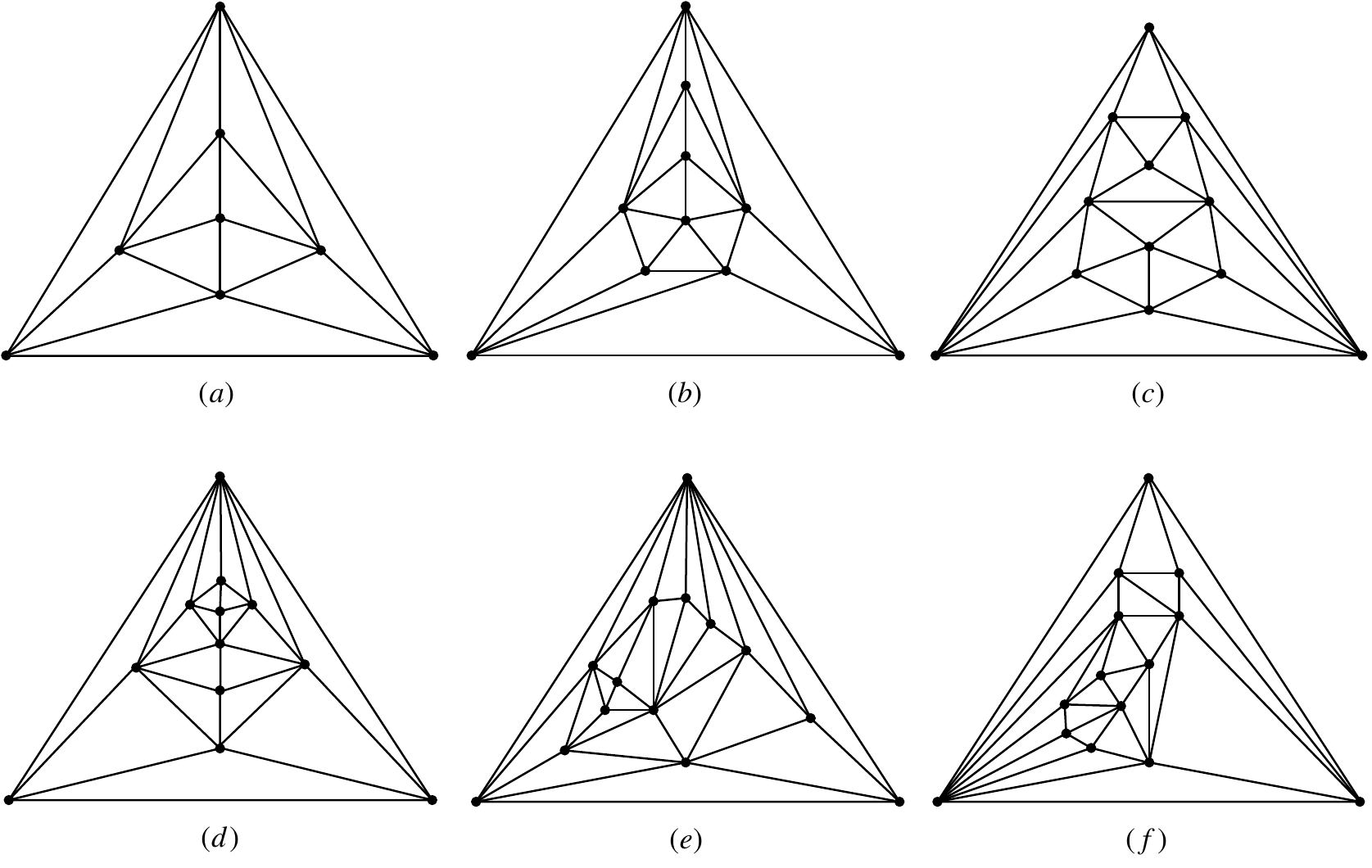}\\
  \caption{Examples of  possible structures of the subgraph induced by the four odd-vertices of a graph in MPG4 }\label{figure2.7}
\end{figure}

Based on the discussion of  Theorem \ref{theorem 2.5} and \ref{theorem2.7}, we have figured out  the impossible structure of the subgraph induced by the four odd-vertices for a graph in $MPG4$. However, all other structures of this subgraph can appear, including a 4-cycle (see Figure \ref{figure2.7} (\emph{a})), a path on 4 vertices  (see Figure \ref{figure2.7} (\emph{b})), two vertex-disjoint $K_2$ (see Figure \ref{figure2.7} (\emph{c})),  a path on 3 vertices and a isolated vertex (see Figure \ref{figure2.7} (\emph{d})), a $K_2$ and two isolated vertices (see Figure \ref{figure2.7} (\emph{e})), and four isolated vertices see Figure \ref{figure2.7} (\emph{f})).

\noindent{\bf{Remark.}}
In this article, we investigated a class of maximal planar graphs, called tree-colorable maximal planar graphs.
We proved that a tree-colorable maximal planar graph $G$ with $\delta(G)\geq 4$ contains at least four odd-vertices.
In addition, for a graph $G$ in $MPG4$, we showed that the subgraph induced by its four odd-vertices is not a claw  and contains no triangles.

With the results we have gained, one can construct maximal planar graphs that contain no tree-colorings.
However, for a given maximal planar graph $G$ that contains exactly four odd-vertices, how to determine whether $G$ is tree-colorable is still unclear. Exploring the sufficient conditions for $G$ to be tree-colorable is an challenging task, which we will research on in the future.


\section*{References}
\begin{thebibliography}{9}


{\small

\bibitem{Bondy2008} J. A. Bondy and U. S. R. Murty, Graph Theory, Springer, 2008.


\bibitem{Borodin1976} O.V. Borodin, A proof of Grunbaum's Conjecture on the acyclic 5-colorability of planar graphs, Dokl. Akad. Nauk SSSR 231(1976), 18-20[Sov. Math. Dokl. 17(1976),1499-1502(1977)].

\bibitem{Borodin1979} O.V. Borodin, On acyclic colorings of planar graphs, Discrete Math, 25(1979), 211-236.

\bibitem{Borodin1999} O.V. Borodin, A.V. Kostochka, and D.R. Woodall, Acyclic colorings of planar graphs with large girth, J. London Math. Soc. 60(1999), 344-352.

\bibitem{Borodin2009} O. V. Borodin, Acyclic 4-colorability of planar graphs without cycles of length 4 or 6, Kiskretn. Anal. Issled. Oper. (6)16(2009), 3-11[J. Appl. Indust. Math. 4(4),490-495(2010)].

\bibitem{Borodin2011} O. V. Borodin, Acyclic 4-coloring of planar graphs without 4- and 5-cycles, Journal of Applied and Industrial Mathematics, (1)5(2011), 31-43.

    \bibitem{Borodin2013} O. V. Borodin, Acyclic 4-Choosability of Planar Graphs with No 4- and 5-Cycles, Journal of Graph Theory, (4)72(2013), 374每397.


 \bibitem{McKay2007}
    G. Brinkmann, B. McKay, Fast generation of planar graphs, MATCH Commum. Math. Comput. Chem., 58(2007), 323-357.

\bibitem{Chen2009} M. Chen and A. Raspaud, Planar graphs without 4, 5 and 8-cycles are acyclically 4-choosable, Electronic Notes in Discrete Mathematics, 34(2009), 659-667.


\bibitem{Garey1976} M. R. Garey, D. S. Johnson and R. Endre Tarjan, The planar Hamiltonian circuit problem is NP-complete, SIAM J. Comput, (4)5(1976),704-714.

\bibitem{Grunbaum1973} B. Gr\"{u}nbaum, Acyclic colorings of planar graphs, Israel J.Math, (3)14(1973), 390-408.

\bibitem{Kostochka1976} A.V. Kostochka and L.S. Melnikov, To the paper of B. Gr\"{u}nbaum on acyclic colorings, Discrete Math. 14(1976), 403-406.


\bibitem{Mondala2012}
    D. Mondal,R. I. Nishat, S. Whitesides, M. S. Rahman, Acyclic colorings of graph subdivisions revisited, Journal of Discrete Algorithms 16(2012), 90每103.

\bibitem{Mondala2013} D. Mondala, R. I. Nishatb, Md.S. Rahmanc, S. Whitesidesb, Acyclic coloring with few division vertices, Journal of Discrete Algorithms, 2013.

\bibitem{Montassier2006} M. Montassier, A. Raspaud, and W. Wang, Acyclic 4-choosability of planar graphs without cycles of specific lengths, In Topics in Discrete Mathematics: Algorithms and Combinatorics, Vol. 26(Springer, Berlin, 2006),pp. 473-491.


 \bibitem{Ochem2005}
    P. Ochem, Negative results on acyclic improper colorings, in:European Conferenceon Combinatorics,EuroComb, in:DMTCSProc.,vol.AE, 2005, pp.357每362.

\bibitem{Wegner1973} G. Wegner, On the paper of B. Gr\"{u}nbaum on  acyclic colorings, Israel J. Math. 14(1973),409-412.




}
\end{thebibliography}
\end{document}